\documentclass[a4paper]{amsart}

\usepackage{amsmath,amssymb,amsthm,url}
\usepackage[utf8]{inputenc}
\usepackage[T1]{fontenc}
\usepackage{libertine}
\usepackage[libertine,cmintegrals,cmbraces]{newtxmath}
\usepackage{verbatim,enumitem}

\AtBeginDocument{%
   \def\MR#1{}
}

\usepackage{stmaryrd}
\usepackage{mathtools}
\usepackage{microtype}
\usepackage{xcolor}
\definecolor{darkred}{RGB}{139,0,0}
\definecolor{darkblue}{RGB}{0,0,139}
\definecolor{darkgreen}{RGB}{0,100,0}
\usepackage{hyperref}
\hypersetup{
  colorlinks   = true, 
  urlcolor     = darkred, 
  linkcolor    = darkred, 
  citecolor   = darkgreen}
\usepackage{tikz}
\usepackage{tikz-cd}
\usepackage[enableskew]{youngtab}
\usepackage{ytableau}
\usepackage[capitalise]{cleveref}
\usepackage{nicefrac}

\tikzset{
  symbol/.style={
    draw=none,
    every to/.append style={
      edge node={node [sloped, allow upside down, auto=false]{$#1$}}}
  }
}

\setenumerate[1]{label=(\roman*),}
\setitemize[1]{label=\raisebox{0.25ex}{\tiny$\bullet$}}

\newcommand{\MTString}{\ensuremath{\mathbf{MTString}}}
\newcommand{\MString}{\ensuremath{\mathbf{MString}}}
\newcommand{\BString}{\ensuremath{\operatorname{BString}}}

\newcommand{\MSO}{\ensuremath{\mathbf{MSO}}}

\newcommand{\HZ}{\ensuremath{\mathbf{HZ}}}

\newcommand{\BO}{\ensuremath{\operatorname{BO}}}

\newcommand{\Diff}{\ensuremath{\operatorname{Diff}}}

\newcommand{\BDiff}{\ensuremath{\operatorname{BDiff}}}

\DeclareMathAlphabet{\mathpzc}{OT1}{pzc}{m}{it}

\newcommand{\oH}{\ensuremath{\operatorname{H}}}

\newcommand{\oO}{\ensuremath{\operatorname{O}}}

\newcommand{\bfR}{\ensuremath{\mathbf{R}}}
\newcommand{\bfZ}{\ensuremath{\mathbf{Z}}}

\newcommand{\bfS}{\ensuremath{\mathbf{S}}}

\newcommand{\ra}{\rightarrow}
\newcommand{\lra}{\longrightarrow}

\newcommand{\xlra}[1]{\overset{#1}{\longrightarrow}}

\newcommand{\interior}[1]{\ensuremath{\operatorname{int}(#1)}}

\newcommand{\K}{\mathrm{K}}
\newcommand{\BGL}{\mathrm{BGL}}
\newcommand{\KSp}{\mathrm{KSp}}
\newcommand{\BSp}{\mathrm{BSp}}
\newcommand{\BGamma}{\mathrm{B\Gamma}}
\newcommand{\Sp}{\mathrm{Sp}}
\newcommand{\GL}{\mathrm{GL}}

\newcommand{\coker}{\mathrm{coker}}
\newcommand{\bP}{\mathrm{bP}}

\newcommand{\sgn}{\mathrm{sgn}}

\newcommand{\Th}{\mathrm{Th}}

\newcommand{\Mod}[1]{\ (\mathrm{mod}\ #1)}

\newtheorem{thm}{Theorem}[section]
\newtheorem*{nthm}{Theorem}
\newtheorem{lem}[thm]{Lemma}

\theoremstyle{definition}

\theoremstyle{remark}

\newtheorem*{nrem}{Remark}

\begin{document}

\title[Some hermitian K-groups via geometric topology]{Some hermitian K-groups via geometric topology}


\author{Manuel Krannich}
\email{krannich@dpmms.cam.ac.uk}
\address{Centre for Mathematical Sciences, Wilberforce Road, Cambridge CB3 0WB, UK}

\author{Alexander Kupers}
\email{kupers@math.harvard.edu}
\address{Department of Mathematics, Harvard University, One Oxford Street, Cambridge MA, 02138, USA}

\begin{abstract}
We compute the first two symplectic quadratic K-theory groups of the integers, or equivalently, the first two stable homology groups of the group of symplectic integral matrices preserving the standard quadratic refinement. The main novelty in our calculation lies in its method, which is based on high-dimensional manifold theory.
\end{abstract}

\maketitle

\section{Introduction}

The \emph{algebraic $K$-theory groups} $\K_*(\bfZ)$ of the integers are the homotopy groups of the algebraic $K$-theory spectrum $\K(\bfZ)$. The zero component of its infinite loop space agrees with the plus construction of the stable general linear group $\BGL_\infty(\bfZ)$, so one has
\[\K_*(\bfZ) \cong  \pi_* \BGL_\infty(\bfZ)^+\quad\text{for }*>0.\]
The \emph{hermitian K-theory groups} of the integers can be defined similarly through the group completion of the automorphism groups of $(\epsilon,\Lambda)$-quadratic modules for a form parameter $(\epsilon,\Lambda)$ à la Bak \cite{Bak}. For instance, in the case $(\epsilon,\Lambda) = (-1,\bfZ)$, the automorphism groups are the symplectic groups $\Sp_{2g}(\bfZ)$, giving rise to the symplectic $K$-theory groups $\KSp_*(\bfZ)$.

Away from the prime $2$, many of these hermitian K-groups have been computed by Karoubi \cite{Karoubi}. More recently, with work of Schlichting \cite{Schlichting} and forthcoming work by Calmès, Dotto, Harpaz, Hebestreit, Land, Moi, Nardin, Nikolaus, and Steimle \cite{9authors}, computations at the prime $2$ became tractable by relating hermitian $K$-theory to the homotopy orbits of a $C_2$-action on $\K(\bfZ)$ and certain $L$-theory spectra.

In this note, we observe that some of these groups can be computed in a very different, somewhat surprising way: via parametrised high-dimensional manifold theory. We carry out the first case that is approachable from this point of view and was not known---the second hermitian $K$-group corresponding to the form parameter $(\epsilon,\Lambda) = (-1,2\bfZ)$. The automorphisms of the associated hyperbolic quadratic modules are the finite index subgroups $\Sp^q_{2g}(\bfZ) \subset \Sp_{2g}(\bfZ)$ of the symplectic groups that preserve the standard quadratic refinement $q \colon \smash{\bfZ^{2g}} \ra \bfZ/2$ of Arf invariant $0$, sending $\sum_{i=1}^g(a_i e_i+b_i f_i)$ to $\sum_{i=1}^ga_ib_i\Mod{2}$ for $a_i,b_i\in\bfZ$ and $(e_i,f_i)_{1\le i\le g}$ the standard symplectic basis of $\bfZ^{2g}$. Hence, the associated hermitian $K$-theory groups are given as
\[\KSp^q_*(\bfZ)\cong\pi_*\BSp^q_\infty(\bfZ)^+\quad\text{for }*>0.\]
The torsion free quotient of $\KSp^q_2(\bfZ)$ is closely related to Meyer's \emph{signature class} $\sgn\in\oH^2(\BSp_{\infty}(\bfZ))$, which measures the signature of total spaces of manifold bundles over surfaces \cite{MeyerThesis, Meyer}. Its pullback to $\BSp_{\infty}^q(\bfZ)$ induces a morphism $\sgn\colon \oH_2(\BSp^q_\infty(\bfZ))\ra\bfZ$ which is not quite an isomorphism, but we show that it is close to it.

\begin{nthm}In degrees $*\le 2$, the Hurewicz map
\[\KSp^q_*(\bfZ)\lra \oH_*(\BSp^q_{\infty}(\bfZ)^+) = \oH_*(\BSp^q_{\infty}(\bfZ))\] is an isomorphism and we have
\[\oH_1(\BSp^q_\infty(\bfZ))\cong\bfZ/4\quad\text{and}\quad\oH_2(\BSp^q_\infty(\bfZ))\cong\bfZ,\] where the second isomorphism is given by the signature morphism divided by $8$.
\end{nthm}

\begin{nrem}\ 
\begin{enumerate}
\item Keeping track of the stability ranges that appear in the proof of the theorem, our argument proves the more refined result that
\[\oH_1(\BSp^q_g(\bfZ))\cong \bfZ/4\text{ for }g\ge 5\quad\text{and}\quad\oH_2(\BSp^q_g(\bfZ))\cong \bfZ\text{ for }g\ge 7.\]
These stability ranges differ from those resulting from a homological stability result of Charney \cite[Cor.\,4.5]{Charney} by one degree, which is due to an improvement of a connectivity range in Charney's argument contained in \cite{GRWI}.
\item The abelianisation $\oH_1(\BSp^q_\infty(\bfZ))$ was computed via different methods in \cite{Endres,JohnsonMillson}, and the group $\oH_2(\BSp^q_\infty(\bfZ))$ can, up to $2$-torsion, be deduced from \cite{Karoubi}. The new part of our theorem is the method of proof, and the lack of $2$-torsion in the second homology group $\oH_2(\BSp^q_\infty(\bfZ))$.
\item For comparison, it is known that \[\KSp_1(\bfZ)\cong \oH_1(\BSp_\infty(\bfZ))=0\quad\text{ and }\quad\KSp_2(\bfZ)\cong\oH_2(\BSp_\infty(\bfZ))\cong\bfZ,\] where the second isomorphism is given by the signature divided by $4$. The first statement follows from the perfectness of $\Sp_\infty(\bfZ)$. The fact that the second group is free of rank $1$ is well-known (see e.g.\,\cite[Thm\,5.1]{Putman}) and the explicit isomorphism in terms of the signature is a consequence of Meyer's work \cite{MeyerThesis,Meyer}.
\end{enumerate}\end{nrem}

\subsection*{Acknowledgements}This note grew out of a question about $\KSp_2^q(\bfZ)$ asked to us by Fabian Hebestreit at the 2019 Clay Research Conference. We would like to thank him for encouraging us to make our computation available and for valuable comments on an earlier version of this note. The first author was supported by O.\,Randal-Williams' Philip Leverhulme Prize from the Leverhulme Trust and by the ERC under the European Union’s Horizon 2020 research and innovation programme (grant agreement No. 756444). The second author was supported by NSF grant DMS-1803766.
\section{The computation}

\subsection{An elementary lemma}\label{sect:lemma}
At several points in the proof, we shall make use of the following direct consequence of the Künneth theorem and the fact that the first $k$-invariant of a double loop space vanishes (see e.g.\cite{ArlettazDouble}). When applied to the infinite loop space $\BSp^q_\infty(\bfZ)^+$, this lemma in particular shows that the first part of the theorem is an immediate consequence of the second, using that the second homology of a cyclic group vanishes.
\begin{lem}\label{lemma}Let $X$ be a connected double loop space. The natural sequence
\[0\lra \pi_2X\xlra{h} \oH_2(X)\xlra{t}\oH_2(K(\pi_1X,1))\lra0\] is split exact, where $h$ is the Hurewicz homomorphism and $t$ is induced by $1$-truncation.
\end{lem}

\subsection{Relation to mapping class groups}
The symplectic group $\smash{\Sp_{2g}}(\bfZ)$ and its finite index \emph{theta subgroup} $\Sp_{2g}^q(\bfZ)\le \Sp_{2g}(\bfZ)$ as defined in the introduction are closely related to the topological group $\Diff_\partial(W^{2n}_{g,1})$ of diffeomorphisms of the manifold \[W^{2n}_{g,1} \coloneqq \sharp^g S^n\times S^n\backslash \interior{D^{2n}}\] for $n$ odd that fix a neighborhood of the boundary pointwise. Indeed, the action of $\Diff_\partial(W^{2n}_{g,1})$ on the homology $\oH_n(W^{2n}_{g,1})\cong\bfZ^{2g}$ gives rise to a homomorphism \[\Diff_\partial(W^{2n}_{g,1})\lra \GL_{2g}(\bfZ)\] whose image agrees with $\Sp_{2g}(\bfZ)$ if $n=1,3,7$ and with $\Sp_{2g}^q(\bfZ)\le \Sp_{2g}(\bfZ)$ for $n\ne1,3,7$ odd. This is well-known for $n=1$ and a consequence of \cite[Prop.\,3]{Kreck} for larger $n$. This morphism evidently factors over the \emph{mapping class group}
\[\Gamma^n_{g,1}\coloneqq \pi_0\Diff_\partial(W^{2n}_{g,1})\] and the kernel of the resulting map $\Gamma^n_{g,1}\ra \GL_{2g}(\bfZ)$ is the \emph{Torelli subgroup of $\Gamma^n_{g,1}$}. In contrast to the case $n=1$, this subgroup is for $n\ge3$ well understood: Kreck \cite{Kreck} determined it up to an extension, and the remaining ambiguity was resolved in \cite{Krannich}. 

In the case $n\equiv 5\Mod{8}$ the Torelli subgroup is particularly simple as it agrees with the subgroup of isotopy classes supported in a codimension zero disc, which is in turn isomorphic to the group $\Theta_{2n+1}$ of homotopy $(2n+1)$-spheres. The resulting extension 
\begin{equation}\label{ses}0\lra \Theta_{2n+1}\lra \Gamma^n_{g,1}\lra \Sp_{2g}^q(\bfZ)\lra 0\end{equation} is central, since every diffeomorphism of $W^{2n}_{g,1}$ fixes a disc up to isotopy.

A priori, the homology of the mapping class group seems no less difficult to compute than that of $\Sp_{2g}^q(\bfZ)$, but work of Galatius--Randal-Williams \cite{GRWstable,GRWI} allows us to access the homology of the classifying space $\BDiff_\partial(W^{2n}_{g,1})$ of the full diffeomorphism group in a range, which we shall use below to deduce information about the homology of $\Gamma^n_{g,1}$ in low degrees and ultimately about that of $\Sp_{2g}^q(\bfZ)$, via the extension \eqref{ses}.

\subsection{Analysing the extension}\label{sect:extension}
The right hand side of \eqref{ses} is independent of the choice of dimension $n\equiv5\Mod{8}$, so to compute something about $\Sp_{2g}^q(\bfZ)$, we may pick whichever dimension of $W^{2n}_{g,1}$ is most convenient. It turns out that $n=5$ is the best choice for our purposes, because certain bordism groups appearing in our calculation are especially well-understood for $n=5$ and the cokernel of the $J$-homomorphism $\coker(J)_{2n+1}$ vanishes in this case, so the group of homotopy spheres $\Theta_{11}$ agrees with its cyclic subgroup $\bP_{12}$ of homotopy spheres that bound parallelizable manifolds (see \cite{KervaireMilnor}). This group is generated by the $11$-dimensional \emph{Milnor sphere} $\Sigma_M\in \bP_{12}$, which is the boundary of the $E_8$-plumbing in dimension $12$ (see \cite[§3]{Levine} and \cite[V.2]{BrowderBook}). In the case $n=5$, the Hochschild--Serre spectral sequence of \eqref{ses} thus induces an exact sequence of the form
\[\oH_2(\BGamma^5_{g,1})\lra \oH_2(\BSp_{2g}^q(\bfZ))\xlra{d_2} \bP_{12}\lra \oH_1(\BGamma^5_{g,1})\lra \oH_1(\BSp_{2g}^q(\bfZ))\lra 0.\] Using the description of the differentials on the $E_2$-page of the Hochschild--Serre spectral sequence of a group extension with abelian kernel in terms of the extension class (c.f.\,\cite[Thm 4]{HochschildSerre}), it follows from \cite[Thm\,B]{Krannich} that the differential $d_2$ is given by 
\[\textstyle{(\sgn/8)\cdot \Sigma_M\colon \oH_2(\BSp_{2g}^q(\bfZ))\lra\bP_{12}},
\] where $\sgn\colon \oH_2(\BSp^q_{2g}(\bfZ))\ra\bfZ$ is the signature morphism (see e.g.\,\cite[Sect.\,3.4]{Krannich}). 

As $\oH_2(\BSp^q_{2g}(\bfZ))$ contains for $g\ge2$ a class with signature $8$ by \cite[Lem.\,3.14]{Krannich}\footnote{See the end of the proof of \cite[Thm\,7.7]{GRWabelian} for a different proof of this fact for $g\ge4$.}, this differential is surjective as long as $g\ge2$. Under this assumption, we thus have an isomorphism $\oH_1(\BGamma^5_{g,1})\cong\oH_1(\BSp_g^{q}(\bfZ))$
and a commutative diagram with exact rows and vertical epimorphisms
\begin{equation}\label{signature}
\begin{tikzcd}
&\oH_2(\BGamma^5_{g,1})\rar\arrow[d,swap, two heads, "\sgn/(8\lvert\bP_{12}\rvert)"]& \oH_2(\BSp_{2g}^q(\bfZ))\arrow[r,"(\sgn/8)\cdot \Sigma_M"]\arrow[d,"\sgn/8", two heads]&[+40pt] \bP_{12}\arrow[d,equal]\arrow[r] &0\\
0\rar& \bfZ\arrow[r,"\lvert\bP_{12}\rvert"]& \bfZ\arrow[r,"{\Sigma_M}"]& \bP_{12}\arrow[r]& 0.
\end{tikzcd}
\end{equation}
Hence, to prove the theorem, it suffices to show that $\oH_1(\BGamma^5_{g,1})\cong\bfZ/4$ and that the left vertical map in the diagram is injective. Instead, we shall consider the full diffeomorphism group $\BDiff_\partial(W^{10}_{g,1})$ and prove $\oH_1(\BDiff_\partial(W^{10}_{g,1})\cong\bfZ/4$ as well as that the divided signature \[\sgn/(8\lvert\bP_{12}\rvert)\colon \oH_2(\BDiff_\partial (W^{10}_{g,1}))\lra \bfZ\] is injective. This implies the previous claims on the homology of $\BGamma^5_{g,1}$ since the canonical map $\BDiff(W^{10}_{g,1})\ra \BGamma^5_{g,1}$ induced by taking path components is $2$-connected.

\subsection{Low-degree homology of $\BDiff_\partial(W^{10}_{g,1})$}The discussion above reduces the proof of the theorem to a question about the homology of $\BDiff_\partial(W^{10}_{g,1})$, which is in a range accessible through the work of Galatius--Randal-Williams \cite{GRWI,GRWstable}. Their results in particular show that the parametrised Pontryagin--Thom collapse map
\begin{equation}\label{paramPT}
\BDiff_\partial(W^{10}_{g,1})\lra \Omega^\infty_0\MTString(10)
\end{equation}
induces an isomorphism in degrees $*\le(g-3)/2$. Here $\MTString(10)$ is the Thom spectrum $\Th(-\theta^*\gamma)$ of the inverse of the pullback of the universal bundle over $\BO(10)$ along the canonical map $\theta\colon\BString(10)\ra\BO(10)$. For $g\ge5$, we therefore obtain an isomorphism $\oH_1(\BDiff_\partial(W^{10}_{g,1}))\cong \pi_1\MTString(10)$ and  for $g\ge7$ a split short exact sequence
 \begin{equation}
 \label{sesMTString}0\lra \pi_2\MTString(10)\lra\oH_2(\BDiff_\partial(W^{10}_{g,1}))\lra\oH_2(B\pi_1\MTString(10))\lra 0,
 \end{equation}
as a consequence of \cref{lemma} applied to $\Omega_0^\infty\MTString(10)$. 

Consequently, to compute the first two homology groups of $\BDiff_\partial(W^{10}_{g,1})$ for large $g$, it suffices to determine the groups $\pi_1\MTString(10)$ and $\pi_2\MTString(10)$. To do so, we adapt a method by Galatius--Randal-Williams \cite[Sect.\,5.1]{GRWabelian}: writing $\MString$ for the Thom spectrum of the universal bundle over $\BString$, there is a canonical stabilisation map \begin{equation}\label{stabilisation}\MTString(10)\lra\Sigma^{-10}\MString\end{equation} whose homotopy fibre admits a $5$-connected map from $\Sigma^{\infty-10}\oO/\oO(10)$. As $\oO/\oO(10)$ is $9$-connected, its unstable and stable homotopy groups agree up to degree $18$ by Freudenthal's suspension theorem, so the long exact sequence induced by \eqref{stabilisation} has  in low degrees the form
\begin{equation}\label{les}
\begin{tikzcd}
\pi_{12}(\oO/\oO(10))\rar& \pi_2\MTString(10)\rar \ar[draw=none]{d}[name=Y, anchor=center]{} &\pi_{12}\MString \ar[rounded corners,
to path={ -- ([xshift=2ex]\tikztostart.east)
	|- (Y.center) \tikztonodes
	-| ([xshift=-2ex]\tikztotarget.west)
	-- (\tikztotarget)}]{dll}\\
\pi_{11}(\oO/\oO(10))\rar& \pi_1\MTString(10)\rar \ar[draw=none]{d}[name=X, anchor=center]{}&\pi_{11}\MString \ar[rounded corners,
to path={ -- ([xshift=2ex]\tikztostart.east)
	|- (X.center) \tikztonodes
	-| ([xshift=-2ex]\tikztotarget.west)
	-- (\tikztotarget)}]{dll}\\
\pi_{10}(\oO/\oO(10))\rar& \pi_0\MTString(10)\rar&\pi_{10}\MString\arrow[r]&0.
\end{tikzcd}
\end{equation}
Comparing $\oO/\oO(10)$ to the Stiefel manifold $\oO(10+k)/\oO(10)$ of $k$-frames in $\bfR^{10+k}$ for large $k$, computations of Paechter \cite{Paechter} (see also \cite{HooMahowald}) show that \[\pi_{12}\oO/\oO(10)=0,\quad \pi_{11}\oO/\oO(10)\cong\bfZ/4,\quad\text{and}\quad\pi_{10}\oO/\oO(10)\cong\bfZ.\] As the group $\pi_{11}\MString$ vanishes by \cite{Giambalvo}, the long exact sequence \eqref{les} yields a short exact sequence of the form
\begin{equation}\label{ses2}
\begin{tikzcd}[column sep=0.4cm]
0\rar& \pi_2\MTString(10)\arrow[r]&\pi_{12}\MString\rar&\bfZ/4\rar& \pi_1\MTString(10)\rar&0.
\end{tikzcd}
\end{equation}

As a next step, we show that the morphism $\pi_2\MTString(10)\ra\pi_{12}\MString$ is surjective. By \cite{Giambalvo}, the group $\pi_{12}\MString$ is free of rank $1$ and we claim that, via the Pontryagin--Thom isomorphism, a generator of this group is represented by a String-manifold of signature $8\lvert\bP_{12}\rvert$. To see this, note that every class in the bordism group $\pi_{12}\MString$ can be represented by almost parallelisable $5$-connected manifold by surgery theory, so its signature is divisible by $8 \lvert\bP_{12}\rvert$ (cf.\,\cite[Ch.\,7]{KervaireMilnor}). On the other hand, the group $\pi_{12}\MString$ does contain such a class of signature $8 \lvert\bP_{12}\rvert$, for instance the class represented by the $12$-manifold obtained from the $\lvert\bP_{12}\rvert$-fold boundary connected sum of the $E_8$-plumbing by gluing a $12$-disc to the boundary. As a result, we have an isomorphism $\pi_{12}\MString\cong\bfZ$ given by the signature divided by $8 \lvert\bP_{12}\rvert$, so we are left with showing that $\pi_2\MTString(10)$ contains a class with that signature. 

To see this, we consider the commutative square
\[
\begin{tikzcd}
\pi_2\MTString(10)\arrow[r]\arrow[d]&\MSO_2(\tau_{>0}\MTString(10))\dar\\
\pi_{12}\MString\arrow[r]&\pi_{12}\MSO.
\end{tikzcd}
\]
In this diagram, the right arrow is the composition of the map
\[\MSO_2(\tau_{>0}\MTString(10)) \lra \MSO_2(\Sigma^{-10} \MSO)\]
induced by the connected cover $\tau_{>0}\MTString(10)\ra \MTString(10)$, the stabilisation map \eqref{stabilisation}, and the canonical map $\MString\ra\MSO$, with the map $\MSO_2(\Sigma^{-10} \MSO)\to \MSO_{12}$ induced by the multiplication of $\MSO$. The upper arrow is induced by the unit $\bfS\ra\MSO$. Since this unit map is $1$-connected and $\tau_{>0} \MTString(10)$ is $0$-connected, the upper arrow is surjective, so it suffices to show that the image of the right vertical map contains a class of signature $8 \lvert\bP_{12}\rvert$. Chasing through the Pontryagin--Thom construction, one sees that the precomposition of the right map with the composition
\[\MSO_2(\BDiff_\partial(W^{10}_{g,1}))\ra\MSO_2(\Omega^{\infty}_0\MTString(10))\ra \MSO_2(\tau_{>0}\MTString(10))\] induced by the counit $\Sigma^\infty\Omega^\infty_0\MTString(10)\simeq \Sigma^\infty\Omega^\infty\tau_{>0}\MTString(10)\ra \tau_{>0}\MTString(10)$ and \eqref{paramPT} maps a class represented by a $W^{10}_{g,1}$-bundle over a surface with trivialised boundary bundle to the (oriented) bordism class of its total space. 

Hence it suffices to show that there is such a bundle whose total space has signature $8\lvert\bP_{12}\rvert$. This follows from the fact that the canonical map $\MSO_2(\BDiff_{\partial}(W^{10}_{g,1}))\ra\oH_2(\BDiff_{\partial}(W^{10}_{g,1}))$ is an isomorphism since $\MSO\ra\HZ$ is $4$-connected, together with the exactness of the first row of \eqref{signature}, because $\oH_2(\BSp_{2g}^q(\bfZ))$ contains a class of signature $8$ and the map $\BDiff_\partial(W^{10}_{g,1})\ra \BGamma^5_{g,1}$ is surjective on second homology as it is $2$-connected.

Having shown that $\pi_2\MTString(10)\ra\pi_{12}\MString$ is surjective, the sequences \eqref{ses2} and \eqref{sesMTString} result in isomorphisms of the form
\begin{align*}\oH_1(\BDiff_\partial(W^{10}_{g,1}))&\cong\pi_1\MTString(10)\cong\bfZ/4\ \ \text{and}\ \\
\oH_2(\BDiff_\partial(W^{10}_{g,1}))&\cong\pi_2\MTString(10)\cong\bfZ,\end{align*}
where the last isomorphism is given by the signature divided by $8\lvert\bP_{12}\rvert$. By the discussion at the end of \cref{sect:extension}, this is sufficient to conclude the theorem.

\begin{nrem}\ 
\begin{enumerate}
\item We remark that our calculation does not rely on the previously known computation of $\oH_1(\BSp_\infty^q(\bfZ))$, even though some arguments in \cite{GRWabelian} and \cite{Krannich} do.
\item For $n$ even, the image of the morphism $\Diff_\partial(W^{2n}_{g,1})\ra\GL_{2g}(\bfZ)$ given by acting on the middle homology agrees with the subgroup $\oO_{g,g}(\bfZ)\subset \GL_{2g}(\bfZ)$ corresponding to the form parameter $(\epsilon,\Lambda)=(1,0)$. Its stable homology in low degrees can be extracted for instance from \cite{GritsenkoHulekSankaran,Behr,HahnOmeara} to be
\[\oH_1(\BO_{\infty,\infty}(\bfZ))\cong\bfZ/2^2\quad\text{and}\quad\oH_2(\BO_{\infty,\infty}(\bfZ))\cong\bfZ/2^3,\] but one might also try to compute these groups via an approach similar to the calculation for $\Sp_\infty^q(\bfZ)$ above. The most promising choice of $n$ in this case is $n=6$, since Kreck's description of $\Gamma^n_{g,1}$ can be used to establish an isomorphism \[\Gamma^n_{g,1}\cong \oO_{g,g}(\bfZ)\times \pi_{13}\MString\]for $n=6$ \cite[Thm\,7.2]{GRWabelian}. Given that $\pi_{13}\MString\cong\bfZ/3$ by \cite{Giambalvo}, it therefore suffices to compute the homology of $\BGamma^6_{g,1}$ in low degrees. Using a similar strategy as in the case $n=5$, one can calculate for large $g$ that \[\oH_1(\BDiff_{\partial}(W^{12}_{g,1})\cong\bfZ/2^2\oplus\bfZ/3\quad\text{and}\quad\oH_2(\BDiff_{\partial}(W^{12}_{g,1})\cong\bfZ/2^{4},\] from which one deduces $\oH_1(\BO_{g,g}(\bfZ))\cong\bfZ/2^2$ and that $\oH_2(\BGamma^6_{g,1})\cong \oH_2(\BO_{g,g}(\bfZ))$ is a quotient of $\bfZ/2^4$. However, narrowing down the isomorphism type of this second homology group further becomes increasingly more difficult than for $\Sp^q_\infty(\bfZ)$ and requires new ideas and inputs, which exceed the scope of this note.\end{enumerate}
\end{nrem}

\bibliographystyle{amsalpha}
\bibliography{literature}

\newcommand{\etalchar}[1]{$^{#1}$}
\providecommand{\bysame}{\leavevmode\hbox to3em{\hrulefill}\thinspace}
\providecommand{\MR}{\relax\ifhmode\unskip\space\fi MR }
\providecommand{\MRhref}[2]{%
  \href{http://www.ams.org/mathscinet-getitem?mr=#1}{#2}
}
\providecommand{\href}[2]{#2}
\begin{thebibliography}{CDH{\etalchar{+}}19}

\bibitem[Arl90]{ArlettazDouble}
D.~Arlettaz, \emph{The first {$k$}-invariant of a double loop space is
  trivial}, Arch. Math. (Basel) \textbf{54} (1990), no.~1, 84--92. \MR{1029601}

\bibitem[Bak81]{Bak}
A.~Bak, \emph{{$K$}-theory of forms}, Annals of Mathematics Studies, vol.~98,
  Princeton University Press, Princeton, N.J.; University of Tokyo Press,
  Tokyo, 1981. \MR{632404}

\bibitem[Beh75]{Behr}
H.~Behr, \emph{Explizite {P}r\"{a}sentation von {C}hevalley-gruppen \"{u}ber
  {Z}}, Math. Z. \textbf{141} (1975), 235--241. \MR{422442}

\bibitem[Bro72]{BrowderBook}
W.~Browder, \emph{Surgery on simply-connected manifolds}, Springer-Verlag, New
  York-Heidelberg, 1972, Ergebnisse der Mathematik und ihrer Grenzgebiete, Band
  65. \MR{0358813}

\bibitem[CDH{\etalchar{+}}19]{9authors}
B.~Calmès, E.~Dotto, Y.~Harpaz, F.~Hebestreit, M.~Land, K.~Moi, D.~Nardin,
  T.~Nikolaus, and W.~Steimle, \emph{{Real K-theory of stable
  infinity-categories}}, in preparation.

\bibitem[Cha87]{Charney}
R.~Charney, \emph{A generalization of a theorem of {V}ogtmann}, Proceedings of
  the {N}orthwestern conference on cohomology of groups ({E}vanston, {I}ll.,
  1985), vol.~44, 1987, pp.~107--125. \MR{885099}

\bibitem[End82]{Endres}
R.~Endres, \emph{Multiplikatorsysteme der symplektischen {T}hetagruppe},
  Monatsh. Math. \textbf{94} (1982), no.~4, 281--297. \MR{685375}

\bibitem[GHS09]{GritsenkoHulekSankaran}
V.~Gritsenko, K.~Hulek, and G.~K. Sankaran, \emph{Abelianisation of orthogonal
  groups and the fundamental group of modular varieties}, J. Algebra
  \textbf{322} (2009), no.~2, 463--478. \MR{2529099}

\bibitem[Gia71]{Giambalvo}
V.~Giambalvo, \emph{On {$\langle 8\rangle$}-cobordism}, Illinois J. Math.
  \textbf{15} (1971), 533--541. \MR{0287553}

\bibitem[GRW14]{GRWstable}
S.~Galatius and O.~Randal-Williams, \emph{Stable moduli spaces of
  high-dimensional manifolds}, Acta Math. \textbf{212} (2014), no.~2, 257--377.
  \MR{3207759}

\bibitem[GRW16]{GRWabelian}
\bysame, \emph{Abelian quotients of mapping class groups of highly connected
  manifolds}, Math. Ann. \textbf{365} (2016), no.~1-2, 857--879. \MR{3498929}

\bibitem[GRW18]{GRWI}
\bysame, \emph{Homological stability for moduli spaces of high dimensional
  manifolds. {I}}, J. Amer. Math. Soc. \textbf{31} (2018), no.~1, 215--264.
  \MR{3718454}

\bibitem[HM65]{HooMahowald}
C.~S. Hoo and M.~E. Mahowald, \emph{Some homotopy groups of {S}tiefel
  manifolds}, Bull. Amer. Math. Soc. \textbf{71} (1965), 661--667. \MR{177412}

\bibitem[HO89]{HahnOmeara}
A.~J. Hahn and O.~T. O'Meara, \emph{The classical groups and {$K$}-theory},
  Grundlehren der Mathematischen Wissenschaften [Fundamental Principles of
  Mathematical Sciences], vol. 291, Springer-Verlag, Berlin, 1989, With a
  foreword by J. Dieudonn\'{e}. \MR{1007302}

\bibitem[HS53]{HochschildSerre}
G.~Hochschild and J-P. Serre, \emph{Cohomology of group extensions}, Trans.
  Amer. Math. Soc. \textbf{74} (1953), no.~1, 110–134.

\bibitem[JM90]{JohnsonMillson}
D.~Johnson and J.~J. Millson, \emph{Modular {L}agrangians and the theta
  multiplier}, Invent. Math. \textbf{100} (1990), no.~1, 143--165. \MR{1037145}

\bibitem[Kar80]{Karoubi}
M.~Karoubi, \emph{Th\'{e}orie de {Q}uillen et homologie du groupe orthogonal},
  Ann. of Math. (2) \textbf{112} (1980), no.~2, 207--257. \MR{592291}

\bibitem[KM63]{KervaireMilnor}
M.~A. Kervaire and J.~W. Milnor, \emph{Groups of homotopy spheres. {I}}, Ann.
  of Math. (2) \textbf{77} (1963), 504--537. \MR{0148075}

\bibitem[Kra19]{Krannich}
M.~Krannich, \emph{{Mapping class groups of highly connected
  $(4k+2)$-manifolds}}, arXiv e-prints (2019).

\bibitem[Kre79]{Kreck}
M.~Kreck, \emph{Isotopy classes of diffeomorphisms of {$(k-1)$}-connected
  almost-parallelizable {$2k$}-manifolds}, {P}roc. {S}ympos., {U}niv. {A}arhus,
  {A}arhus, 1978, Lecture Notes in Math., vol. 763, Springer, Berlin, 1979,
  pp.~643--663. \MR{561244}

\bibitem[Lev85]{Levine}
J.~P. Levine, \emph{Lectures on groups of homotopy spheres}, Algebraic and
  geometric topology ({N}ew {B}runswick, {N}.{J}., 1983), Lecture Notes in
  Math., vol. 1126, Springer, Berlin, 1985, pp.~62--95. \MR{802786}

\bibitem[Mey72]{MeyerThesis}
W.~Meyer, \emph{Die {S}ignatur von lokalen {K}oeffizientensystemen und
  {F}aserb\"undeln}, Bonn. Math. Schr. (1972), no.~53, viii+59. \MR{0305402}

\bibitem[Mey73]{Meyer}
\bysame, \emph{Die {S}ignatur von {F}l\"{a}chenb\"{u}ndeln}, Math. Ann.
  \textbf{201} (1973), 239--264. \MR{331382}

\bibitem[Pae56]{Paechter}
G.~F. Paechter, \emph{The groups {$\pi _{r}(V_{n,\,m})$}. {I}}, Quart. J. Math.
  Oxford Ser. (2) \textbf{7} (1956), 249--268. \MR{0131878}

\bibitem[Put12]{Putman}
A.~Putman, \emph{The {P}icard group of the moduli space of curves with level
  structures}, Duke Math. J. \textbf{161} (2012), no.~4, 623--674. \MR{2891531}

\bibitem[Sch19]{Schlichting}
M.~Schlichting, \emph{Symplectic and orthogonal {$K$}-groups of the integers},
  C. R. Math. Acad. Sci. Paris \textbf{357} (2019), no.~8, 686--690.
  \MR{4015334}

\end{thebibliography}

\vspace{.5cm}

\end{document}